# On R Degrees of Vertices and R Indices of Graphs

Süleyman Ediz

suleymanediz@yyu.edu.tr

Faculty of Education, Yüzüncü Yıl University, Van, Turkey

**Abstract**

Topological indices have been used to modeling biological and chemical properties of molecules in quantitive structure property relationship studies and quantitive structure activity studies. All the degree based topological indices have been defined via classical degree concept. In this paper we define a novel degree concept for a vertex of a simple connected graph: R degree. And also we define R indices of a simple connected graph by using the R degree concept. We compute the R indices for well-known simple connected graphs such as paths, stars, complete graphs and cycles.

**Keywords:** R degree; R indices; Topological indices; QSAR; QSPR

1. Introduction

Graph theory has many applications to modeling real world situations from the basic sciences to engineering and social sciences. Chemical graph theory has an important effect on the development of the chemical sciences by using topological indices. A topological index, which is a graph invariant it does not depend on the labeling or pictorial representation of the graph, is a numerical parameter mathematically derived from the graph structure. The topological indices of molecular graphs are widely used for establishing correlations between the structure of a molecular compound and its physicochemical properties or biological activity. These indices are used in quantitive structure property relations (QSPR) research. Topological indices are important tools for analyzing some physicochemical properties of molecules without performing any experiment. The first distance based topological index was proposed by (Wiener 1947) for modeling physical properties of alcanes, and after him, hundred topological indices were defined by chemists and mathematicians and so many properties of chemical structures were studied. More than forty years ago (Gutman & Trinajstić 1971) defined Zagreb indices which are degree based topological indices. These topological indices were proposed to be measures of branching of the carbon-atom skeleton in (Gutman et al. 1975). The Randic and Zagreb indices are the most used topological indices in chemical and mathematical literature so far. For detailed discussions of both these indices and other well-known topological indices, we refer the interested reader [4-14] and [33-38] and references therein. In 1998, Estrada *et al* modelled the enthalpy of formation of alkanes by using atom- bond connectivity(ABC) index [15]. The ABC index for a connected graph $G$ defined as;

$$ABC(G) = \sum_{uv \in E(G)} \sqrt{\frac{\deg(u)+\deg(v)-2}{\deg(u) \times \deg(v)}} \qquad (1)$$

There are many open problems related to ABC index in the mathematical chemistry literature. We refer the interested reader the sudies of the last two years [16-20].

In 2009, Vukičević and Furtula defined geometric-arithmetic(GA) index and compared GA index with the well-known Randić index [21]. The authors showed that the GA index give better correlation to modelling standard enthalpy of vaporization of octane isomers. The GA index for a connected graph $G$ defined as;

$$GA(G) = \sum_{uv \in E(G)} \frac{2\sqrt{\deg(u) \times \deg(v)}}{\deg(u)+\deg(v)} \qquad (2)$$

After that many studies related to GA index were conducted in view of mathematical chemistry and QSPR researches [22-26].

The harmonic index was defined by Zhong in 2012 [27]. The H index for a connected graph $G$ defined as;

$$H(G) = \sum_{uv \in E(G)} \frac{2}{\deg(u)+\deg(v)} \qquad (3)$$

The relationships between harmonic index and domination like parameters were investigated by Li *et al.* [28] We refer the interested reader for the article related to H index by Ilić and the references therein [29].

The sum-connectivity index ($\chi$) were defined by Zhou and Trinajstić in 2009 [30]. The $\chi$ index for a connected graph $G$ defined as;

$$\chi(G) = \sum_{uv \in E(G)} (\deg(u) + \deg(v))^{-1/2} \qquad (4)$$

(Farahani 2015) computed the sum-connectivity index of carpa designed cycle [31] and (Akhter et al. 2016) investigated the sum-connectivity index of cacti [32].

As of now in the chemical and mathematical literature all degree based topological indices have been defined by using classical degree concept.

In this study our aim is to define a novel degree concept namely, R degrees. Also by using R degrees, we define the first, the second and the third R indices for a simple connected graph.

**2. R Degrees and R Indices**

In this section we give basic definitions and facts about above mentioned graph invariants. A graph $G = (V, E)$ consists of two nonempty sets $V$ and 2-element subsets of $V$ namely $E$. The elements of $V$ are called vertices and the elements of $E$ are called edges. For a vertex $v$, $\deg(v)$ show the number of edges that incident to $v$. The set of all vertices which adjacent to $v$ is called the open neighborhood of $v$ and denoted by $N(v)$. If we add the vertex $v$ to $N(v)$, then we get the closed neighborhood of $v$, $N[v]$.

For a vertex $v$, $S_v = \sum_{u \in N(v)} deg(u)$. For conveince, we name $S_v$ as "the sum degree of $v$" or briefly "sum degree". For a vertex $v$, $M_v = \prod_{u \in N(v)} deg(u)$. For conveince, we name $M_v$ as "the multiplication degree of $v$" or briefly "multiplication degree".

**Definition 2.1.** The R degree of a vertex $v$ of a simple connected graph $G$ defined as;

$$r(v) = M_v + S_v \tag{5}$$

**Definition 2.3.** The first R index of a simple connected graph $G$ defined as;

$$R^1(G) = \sum_{v \in V(G)} r(v)^2 \tag{6}$$

**Definition 2.4.** The second S index of a simple connected graph $G$ defined as;

$$R^2(G) = \sum_{uv \in E(G)} r(u)r(v) \tag{7}$$

**Definition 2.5.** The third S index of a simple connected graph $G$ defined as;

$$R^3(G) = \sum_{uv \in E(G)} [r(u) + r(v)] \tag{8}$$

**Proposition 2.6.** Let $K_n$ be a complete graph with $n$ vertices ($n \geq 3$). Then;

a. $R^1(K_n) = n.\left((n-1)^2((n-1)^{n-3} + 1)\right)^2$

b. $R^2(K_n) = \frac{1}{2}n(n-1)^5((n-1)^{n-3} + 1)^2$

c. $R^3(K_n) = n(n-1)^3((n-1)^{n-3} + 1)$

**Proof.** Let $v \in K_n$. Then $S_v = (n-1)(n-1) = (n-1)^2$ and $M_v = (n-1)^{n-1}$. Therefore $r(v) = |M_v + S_v| = (n-1)^2((n-1)^{n-3} + 1)$. We can begin to compute since all the vertices of complete graph have same R degree.

a. $R^1(K_n) = \sum_{v \in V(K_n)} r(v)^2 = n.\left((n-1)^2((n-1)^{n-3} + 1)\right)^2$

b. $R^2(K_n) = \sum_{uv \in E(K_n)} r(u)r(v) = \frac{n(n-1)}{2}\left((n-1)^2((n-1)^{n-3}+1)\right)^2$

$$= \frac{1}{2}n(n-1)^5((n-1)^{n-3}+1)^2$$

c. $R^3(K_n) = \sum_{uv \in E(K_n)} [r(u)+r(v)] = \frac{n(n-1)}{2} 2(n-1)^2((n-1)^{n-3}+1)$

$$= n(n-1)^3((n-1)^{n-3}+1)$$

**Proposition 2.7.** Let $C_n$ be a cycle graph with $n$ vertices ($n \geq 3$). Then;

a. $R^1(C_n) = 64n$

b. $R^2(C_n) = 64n$

c. $R^3(C_n) = 16n$

**Proof.** Let $v \in C_n$. Then $S_v = 2 + 2 = 4$ and $M_v = 2.2 = 4$. Therefore $r(v) = |M_v + S_v| = 8$ and . We can begin to compute since all the vertices of cycle graph have same R degrees.

a. $R^1(C_n) = \sum_{v \in V(C_n)} r(v)^2 = 64n$

b. $R^2(C_n) = \sum_{uv \in E(C_n)} r(u)r(v) = 64n$

c. $R^3(C_n) = \sum_{uv \in E(C_n)} [r(u)+r(v)] = 16n$

**Proposition 2.8.** Let $P_n$ be a path graph with $n$ vertices ($n \geq 3$). Then;

a. $R^1(P_n) = n + 5/2$

b. $R^2(P_n) = n + 1$

c. $R^3(P_n) = 2n - 10/3$

**Proof.** Let $V(P_n) = \{v_1, v_2, \dots, v_n\}$ and $E(P_n) = \{v_1v_2, v_2v_3, \dots, v_{n-2}v_{n-1}, v_{n-1}v_n\}$. Then $S_{v_1} = S_{v_n} = 2$ and $M_v = 2$. Therefore $r(v_1) = r(v_n) = 4$. Also $S_{v_2} = S_{v_{n-1}} = 3$ and $M_{v_2} = M_{v_{n-1}} = 2$

. Therefore $r(v_1) = r(v_n) = 5$ and $S_{v_3} = S_{v_4} = \cdots = S_{v_{n-3}} = S_{v_{n-2}} = 4$ and $M_{v_3} = M_{v_4} = \cdots = M_{v_{n-3}} = M_{v_{n-2}} = 4$. Therefore $r(v_3) = r(v_4) = \cdots = r(v_{n-3}) = r(v_{n-2}) = 8$. Note that $P_n$ has $n$-vertex and $n-1$ edges. We can begin our computations.

a. $R^1(P_n) = \sum_{v \in V(P_n)} r(v)^2 = 2.5^2 + (n-2).64 = 64n - 78$

b. $R^2(P_n) = \sum_{uv \in E(P_n)} r(u)r(v) = 2.5.8 + (n-3).64 = 64n - 112$

c. $R^3(P_n) = \sum_{uv \in E(P_n)} [r(u) + r(v)] = 2.13 + (n-3).16 = 16n - 22$

**Proposition 2.8.** Let $S_n$ be a star graph with $n$ vertices ($n \geq 3$). Then;

a. $R^1(S_n) = n^2$

b. $R^2(S_n) = 2n(n-1)^2$

c. $R^3(S_n) = (n-1)(3n-2)$

**Proof.** There are two kinds of vertices. For the central vertex $v$, $\deg(v) = n - 1$ and for any pendent vertex $u$, $\deg(u) = 1$. Then $S_v = n - 1$ and $M_v = 1$. Also for any pendent vertex $u$, $S_u = n - 1$ and $M_u = n - 1$. Therefore $r(v) = n$. Also $r(u) = 2n - 2$.

Note that $S_n$ has $n$-vertex and $n - 1$ edges. We can begin our computations.

a. $R^1(S_n) = \sum_{v \in V(S_n)} r(v)^2 = n^2$

b. $R^2(S_n) = \sum_{uv \in E(S_n)} r(u)r(v) = 2n(n-1)^2$

c. $R^3(S_n) = \sum_{uv \in E(S_n)} [r(u) + r(v)] = (n-1)(3n-2)$

**3.Conclusion**

There are many problems for further study about the R indices. The mathematical properties and relations between R indices and other topological indices are interesting problems worth to study. Also, QSPR analysis of R indices may be attract the attention of some mathematical chemists.